\documentclass[12pt]{article}

\begin{document}

\begin{center}{\large \bf  On a Boundary-Value Problem for One Class of Differential Equations of the Fourth Order with Operator Coefficients}
\end{center}

\begin{center}{\bf A.R. Aliev}
\\ {\it Institute of Mathematics and Mechanics of NAS of Azerbaijan,}
\\ {\it Baku State University}
\\ {\it E-mail: alievaraz@yahoo.com}
\end{center}

\begin{center}{\bf 1. Introduction}
\end{center}

An interest to investigations of the initial and boundary-value problems for the operator-differential equations (see, for example, [1-3] and given references there) has increased last years. This is connected with the fact that the equations of this type allow to consider both ordinary differential operators and the operators with partial derivatives.

In this paper the boundary-value problem for one class operator-differential equations of the fourth order is considered in the weighted analogue of Sobolev type space (with the weight $e^{-\frac{\kappa }{2} t} $, $\kappa \in R=(-\infty ;+\infty )$). Namely, the correctness and unique solvability of the boundary-value problem on semi-axis for the operator-differential equation of the fourth order, the main part of which has a multiple characteristic, is studied here. The obtained solvability conditions are expressed in terms of the operator coefficients of the equation, and this allows to check these conditions both in theoretical and in applied problems. Estimations of the norms of the operators of intermediate derivatives closely connected with the solvability conditions have been carried out. Moreover, the connection between the weight exponent and the lower boundary of spectrum of the operator, participating in the equation, is determined in results of the paper.

Let $A$ be a selfadjoint positively defined operator in separable Hilbert space $H$ and $\kappa $ be any real number: $\kappa \in R$.

We denote by $L_{2,\kappa } (R;H)$ Hilbert space of $H$--valued functions, defined in $R$, with the norm

\[\, \left\| f\right\| _{L_{2,\kappa } (R;H)} =\left(\int _{-\infty }^{+\infty }\left\| f(t)\right\| _{H}^{2} e^{-\kappa t} dt \right)^{{\raise0.7ex\hbox{$ 1 $}\!\mathord{\left/{\vphantom{1 2}}\right.\kern-\nulldelimiterspace}\!\lower0.7ex\hbox{$ 2 $}} } .\]
We denote by $W_{2,\kappa }^{4} (R;H)$ the space of $H$--valued functions such that $\frac{d^{4} u(t)}{dt^{4} } \in L_{2,\kappa } (R;H)$, $A^{4} \, u(t)\in L_{2,\kappa } (R;H)$ with the norm

\[\left\| u\right\| _{W_{2,\kappa }^{4} (R;H)} =\left(\left\| \frac{d^{4} u}{dt^{4} } \right\| _{L_{2,\kappa } (R;H)}^{2} +\left\| A^{4} u\right\| _{L_{2,\kappa } (R;H)}^{2} \right)^{{\raise0.7ex\hbox{$ 1 $}\!\mathord{\left/{\vphantom{1 2}}\right.\kern-\nulldelimiterspace}\!\lower0.7ex\hbox{$ 2 $}} } .\]
Here and further the derivatives are considered in sense of the theory of generalized functions. It is obvious that for $\kappa =0$ we'll have the spaces $L_{2,0} (R;H)=L_{2} (R;H)$, $W_{2,0}^{4} (R;H)=W_{2}^{4} (R;H)$ (see [4]). By the same way we can define the spaces $L_{2,\kappa } (R_{+} ;H)$ and $W_{2,\kappa }^{4} (R_{+} ;H)$, where $R_{+} =\left[0;+\infty \right)$.

Let's pass to the statement of the investigated problem. We consider the following boundary-value problem in the space $H$:

\begin{equation} \label{GrindEQ__1_}
\left(-\frac{d}{dt} +A\right)\left(\frac{d}{dt} +A\right)^{3} u(t)+\sum _{j=1}^{4}A_{j} \frac{d^{4-j} u(t)}{dt^{4-j} }  =f(t),\, \, \, t\in R_{+} ,
\end{equation}

\begin{equation} \label{GrindEQ__2_}
u\left(0\right)=\frac{du\left(0\right)}{dt} =\frac{d^{2} u\left(0\right)}{dt^{2} } =0,
\end{equation}
where $A$ is the same operator with the lower boundary of the spectrum $\lambda _{0} $ ($A=A^{*} \ge \lambda _{0} E$ $\, (\lambda _{0} >0)$, $E$ is the identity operator), $A_{j} $, $j=1,2,3,4$ are the linear, generally speaking, unbounded operators, $f(t)\in L_{2,\kappa } (R_{+} ;H)$, $u(t)\in W_{2,\kappa }^{4} (R_{+} ;H)$.

\textbf{Definition.} If for any $f(t)\in L_{2,\kappa } (R_{+} ;H)$ there exists the vector-function $u(t)\in W_{2,\kappa }^{4} (R_{+} ;H)$, satisfying the equation (1) almost everywhere in $R_{+} $, and the boundary conditions (2) in sense:

\[\mathop{\lim }\limits_{t\to 0} \left\| A^{{\raise0.5ex\hbox{$\scriptstyle 7 $}\kern-0.1em/\kern-0.15em\lower0.25ex\hbox{$\scriptstyle 2 $}} -j} \frac{d^{j} u\left(t\right)}{dt^{j} } \right\| _{H} =0,\, \, \, j=0,1,2,\]
and the inequality

\[\left\| u\right\| _{W_{2,\kappa }^{4} (R_{+} ;H)} \le const\left\| f\right\| _{L_{2,\kappa } (R_{+} ;H)} \]
takes place, then it is called the \textit{regular solution} of the boundary-value problem (1), (2), and  the boundary-value problem (1), (2)
is called \textit{regularly solvable}.

We note that similar problems of the regular solvability in the weighted space for the operator-differential equations of the second and third orders, considered on semi-axis, are studied, for example, in the papers [5-8]. Investigations of the regular solvability of the boundary-value problem (1), (2) for $A_{4} =0$ have been done in wide aspect for the case $\kappa =0$ in the paper [9]. But in the paper [10] the sufficient conditions of the normal solvability of the boundary-value problem for differed from our case another class of the operator-differential equations of the fourth order with the multiple characteristic have been obtained in the weighted space. We note that in all these works the solvability conditions are expressed in terms of the operator coefficients of studied equations. In the case of the weighted space it is important to note the earlier paper [11], in which the operators of the perturbed part of the equation are the degrees of the operator $A$, multiplied by the complex numbers, and solvability conditions are expressed with the help of limitations of resolvent increase of the corresponding operator pencil. The solvability and Fredholm solvability of boundary value problems on the semi-axis (as well as on a finite interval) for equations of arbitrary order in Hilbert space with non-commuting operator coefficients was investigated in papers [12, 13]. We also remark that equations of form (1) appear in applications, in particular, in the problems of stability of the plates from the plastic material (see [14]).

\begin{center}{\bf 2. Main results}
\end{center}

We'll begin from studying the operator $P_{0} $, acting from the space $\mathop{W_{2,\kappa }^{4} }\limits^{o} (R_{+} ;H)$ to the space $L_{2,\kappa } (R_{+} ;H)$ by the following way:

\[P_{0} u\left(t\right)\equiv \left(-\frac{d}{dt} +A\right)\left(\frac{d}{dt} +A\right)^{3} u(t),\, \, u\left(t\right)\in \mathop{W_{2,\kappa }^{4} }\limits^{o} \left(R_{+} ;H\right),\]
where

\[\mathop{W_{2,\kappa }^{4} }\limits^{o} (R_{+} ;H)=\left\{u\left(t\right):\, u\left(t\right)\in W_{2,\kappa }^{4} (R_{+} ;H),\, \, \frac{d^{s} u\left(0\right)}{dt^{s} } =0,\, \, s=0,1,2\right\}.\]

It takes place the following theorem on the isomorphism of the operator $P_{0} $, in the proof of which we apply Fourier transform and Banach theorem on the inverse operator.

\textbf{Theorem 1.} \textit{Let $\left|\kappa \right|<2\lambda _{0} $. Then the operator $P_{0} $ isomorphically maps the space $\mathop{W_{2,\kappa }^{4} }\limits^{o} (R_{+} ;H)$ onto the space $L_{2,\kappa } (R_{+} ;H)$.}

\textbf{Proof.} For convenience of the further notes we consider the polynomial operator pencil

\[P_{0} \left(\mu ;A\right)=\left(-\mu E+A\right)\left(\mu E+A\right)^{3} .\]
Then the boundary-value problem (1), (2) for $A_j = 0$, $j = 1,2,3,4$ can be written in the form of the operator equation

\begin{equation} \label{GrindEQ__3_}
P_{0} \left(\frac{d}{dt} ;A\right)u\left(t\right)=f\left(t\right),
\end{equation}
where $f(t)\in L_{2,\kappa } (R_{+} ;H)$, $u(t)\in \mathop{W_{2,\kappa }^{4} }\limits^{o} (R_{+} ;H)$. It is easy to determine that the homogeneous equation $P_{0} \left(\frac{d}{dt} ;A\right)u\left(t\right)=0$ has only trivial solution from the space $\mathop{W_{2,\kappa }^{4} }\limits^{o} (R_{+} ;H)$. Let's show that the equation (3) has the solution from the space $\mathop{W_{2,\kappa }^{4} }\limits^{o} (R_{+} ;H)$ for any $f(t)\in L_{2,\kappa } (R_{+} ;H)$. After substitution $v\left(t\right)=u\left(t\right)e^{-\frac{\kappa }{2} t} $ we rewrite the equation (3) in the form

\begin{equation} \label{GrindEQ__4_}
P_{0} \left(\frac{d}{dt} +\frac{\kappa }{2} ;A\right)v\left(t\right)=g\left(t\right),
\end{equation}
where $v(t)\in \mathop{W_{2}^{4} }\limits^{o} (R_{+} ;H)$, $g\left(t\right)=f(t)e^{-\frac{\kappa }{2} t} \in L_{2} (R_{+} ;H)$. Let $\lambda \in \sigma \left(A\right)$ ($\lambda \ge \lambda _{0} $). As for $\left|\kappa \right|<2\lambda _{0} $

\[\left|P_{0} \left(i\xi +\frac{\kappa }{2} ;\lambda \right)\right|=\left|\left(-\left(i\xi +\frac{\kappa }{2} \right)+\lambda \right)\left(i\xi +\frac{\kappa }{2} +\lambda \right)^{3} \right|=\]

\[\left|\left(-\left(i\xi +\frac{\kappa }{2} \right)^{2} +\lambda ^{2} \right)\left(i\xi +\frac{\kappa }{2} +\lambda \right)^{2} \right|=\]

\[\left(\left(\xi ^{2} -\frac{\kappa ^{2} }{4} +\lambda ^{2} \right)^{2} +\xi ^{2} \kappa ^{2} \right)^{{\raise0.7ex\hbox{$ 1 $}\!\mathord{\left/{\vphantom{1 2}}\right.\kern-\nulldelimiterspace}\!\lower0.7ex\hbox{$ 2 $}} } \left(\left(\lambda +\frac{\kappa }{2} \right)^{2} +\xi ^{2} \right)\ge \]

\[\left(\lambda ^{2} -\frac{\kappa ^{2} }{4} \right)\left(\lambda +\frac{\kappa }{2} \right)^{2} \ge \left(\lambda _{0}^{2} -\frac{\kappa ^{2} }{4} \right)\left(\lambda _{0} +\frac{\kappa }{2} \right)^{2} >0,\,\,\, \xi \in R,\]
then from the spectral decomposition of the operator $A$ it follows that the operator pencil $P_{0} \left(i\xi +\frac{\kappa }{2} ;A\right)$ is invertible for $\left|\kappa \right|<2\lambda _{0} $.

Let us continue the function $f(t)$ by 0 for $t<0$, then (4) will be already on the whole axis and also $g(t)=0$ for $t<0$. Using direct and inverse Fourier transforms it becomes clear that

\[v_{0} \left(t\right)=\frac{1}{2\pi } \int _{-\infty }^{+\infty }P_{0}^{-1} \left(i\xi +\frac{\kappa }{2} ;A\right) \left(\int _{0}^{+\infty }g\left(s\right)e^{-i\xi s} ds \right)e^{i\xi t} d\xi ,\,\,\,  t\in R\]
satisfies the equation (4) almost everywhere in $R$. We'll prove that $v_{0} (t)\in W_{2}^{4} (R;H)$.

Really, from Plancherel theorem we have

\[\left\| v_{0} \right\| _{W_{2}^{4} (R;H)}^{2} =\left\| \frac{d^{4} v_{0} }{dt^{4} } \right\| _{L_{2} (R;H)}^{2} +\left\| A^{4} v_{0} \right\| _{L_{2} (R;H)}^{2} =\left\| \xi ^{4} \hat{v}_{0} \left(\xi \right)\right\| _{L_{2} (R;H)}^{2} +\]

\[\left\| A^{4} \hat{v}_{0} \left(\xi \right)\right\| _{L_{2} (R;H)}^{2} \le \mathop {\sup }\limits_{\xi  \in R} \left\| \xi ^{4} P_{0}^{-1} \left({i\xi
 +\frac{\kappa }{2} ;A}
\right)\right\| _{H\to H}^{2} \left\| \hat{g}\left(\xi
\right)\right\| _{L_{2} (R;H)}^{2} +\]

\[\mathop {\sup }\limits_{\xi  \in R}\left\| A^{4} P_{0}^{-1} \left({i\xi +\frac{\kappa }{2} ;A}\right)\right\| _{H\to H}^{2} \left\| \hat{g}\left(\xi \right)\right\| _{L_{2} (R;H)}^{2} \le \]

\[const\left\| \hat{g}\left(\xi \right)\right\| _{L_{2} (R;H)}^{2} =const\left\| g\left(t\right)\right\| _{L_{2} (R_{+} ;H)}^{2} ,\]
where $\hat{v}_{0} \left(\xi \right)$ and $\hat{g}\left(\xi \right)$
are Fourier transforms of the functions $v_{0} \left(t\right)$ and
$g\left(t\right)$ correspon-dingly. And this in turn so as
estimating the norm $\left\| A^{4} P_{0}^{-1} \left({i\xi
+\frac{\kappa }{2} ;A}\right)\right\|$ for $\xi \in R$, from the
spectral theory of selfadjoint operators, we have

\[\left\| {A^4 P_0^{ - 1} \left( {i\xi  + \frac{\kappa }{2};A} \right)} \right\| = \mathop {\sup }\limits_{\lambda  \in \sigma \left( A \right)} \left| {\lambda ^4 \left( { - \left( {i\xi  + \frac{\kappa }{2}} \right) + \lambda } \right)^{ - 1} \left( {i\xi  + \frac{\kappa }{2} + \lambda } \right)^{ - 3} } \right| =\]

\[\mathop {\sup }\limits_{\lambda  \in \sigma \left( A \right)} \left| {\lambda ^4 \left( { - \left( {i\xi  + \frac{\kappa }{2}} \right)^2  + \lambda ^2 } \right)^{ - 1} \left( {i\xi  + \frac{\kappa }{2} + \lambda } \right)^{ - 2} } \right| \le\]

\[\mathop {\sup }\limits_{\lambda  \in \sigma \left( A \right)} \frac{{\lambda ^4 }}{{\left( {\xi ^2  + \lambda ^2  - \frac{{\kappa ^2 }}{4}} \right)\left( {\left( {\lambda  + \frac{\kappa }{2}} \right)^2  + \xi ^2 } \right)}} \le \mathop {\sup }\limits_{\lambda  \in \sigma \left( A \right)} \frac{{\lambda ^4 }}{{\left( {\lambda ^2  - \frac{{\kappa ^2 }}{4}} \right)\left( {\lambda  + \frac{\kappa }{2}} \right)^2 }} \le\]

\[\frac{{\lambda _0^4 }}{{\left( {\lambda _0^2  - \frac{{\kappa ^2 }}{4}} \right)\left( {\lambda _0  + \frac{\kappa }{2}} \right)^2 }}.\]
Similarly, we prove that for $\xi \in R$

\[\left\| {\xi ^4 P_0^{ - 1} \left( {i\xi  + \frac{\kappa }{2};A} \right)} \right\| = \mathop {\sup }\limits_{\lambda  \in \sigma \left( A \right)} \left| {\xi ^4 \left( { - \left( {i\xi  + \frac{\kappa }{2}} \right) + \lambda } \right)^{ - 1} \left( {i\xi  + \frac{\kappa }{2} + \lambda } \right)^{ - 3} } \right| \leq\]

\[\mathop {\sup }\limits_{\lambda  \in \sigma \left( A \right)} \frac{{\xi ^4 }}{{\left( {\xi ^2  + \lambda ^2  - \frac{{\kappa ^2 }}{4}} \right)\left( {\left( {\lambda  + \frac{\kappa }{2}} \right)^2  + \xi ^2 } \right)}} \le \frac{{\xi ^4 }}{{\left( {\xi ^2  + \lambda _0^2  - \frac{{\kappa ^2 }}{4}} \right)\left( {\left( {\lambda _0  + \frac{\kappa }{2}} \right)^2  + \xi ^2 } \right)}} \le 1.\]
Here ${\sigma \left( A \right)}$ is the spectrum of operator $A$.

As the mapping $v\left(t\right)\to u\left(t\right)e^{-\frac{\kappa }{2} t} $ is the isomorphism between the spaces $W_{2}^{4} (R;H)$ and $W_{2,\kappa }^{4} (R;H)$, then it is obvious that for any $f(t)\in L_{2,\kappa } (R;H)$ there exists $u_{0} \left(t\right)=v_{0} \left(t\right)e^{\frac{\kappa }{2} t} \in W_{2,\kappa }^{4} \left(R;H\right)$, satisfying the equation (3) almost everywhere in $R$. Continuing, we denote by $\tilde{u}_{0} \left(t\right)$ the restriction of a vector-function $u_{0} \left(t\right)$ on $R_{+} $. As $\tilde{u}_{0} \left(t\right)\in W_{2,\kappa }^{4} \left(R_{+} ;H\right)$, then from theorem on the traces [4, chapter 1] $\frac{d^{s} \tilde{u}_{0} \left(0\right)}{dt^{s} } \in D\left(A^{{\raise0.5ex\hbox{$\scriptstyle 7 $}\kern-0.1em/\kern-0.15em\lower0.25ex\hbox{$\scriptstyle 2 $}} -s} \right)$, $s=0,1,2$. As we look for the solution of the boundary-value problem (1), (2) for $A_j = 0$, $j = 1,2,3,4$ in the form

\[u\left(t\right)=\tilde{u}_{0} \left(t\right)+e^{-tA} \varphi _{0} +tAe^{-tA} \varphi _{1} +t^{2} A^{2} e^{-tA} \varphi _{2} ,\]
where $\varphi _{0} ,\varphi _{1} ,\varphi _{2} \in D\left(A^{{\raise0.5ex\hbox{$\scriptstyle 7 $}\kern-0.1em/\kern-0.15em\lower0.25ex\hbox{$\scriptstyle 2 $}} } \right)$, and $e^{-tA} $ is strongly continuous semigroup of the bounded operators, generated by the operator $-A$, then from the boundary conditions (2) we have

\[\left\{\begin{array}{c} {u\left(0\right)=\tilde{u}_{0} \left(0\right)+\varphi _{0} =0,} \\ {\frac{du\left(0\right)}{dt} =\frac{d\tilde{u}_{0} \left(0\right)}{dt} -A\varphi _{0} +A\varphi _{1} =0,} \\ {\frac{d^{2} u\left(0\right)}{dt^{2} } =\frac{d^{2} \tilde{u}_{0} \left(0\right)}{dt^{2} } +A^{2} \varphi _{0} -2A^{2} \varphi _{1} +2A^{2} \varphi _{2} =0,} \end{array}\right. \]
and from here we obtain

\begin{equation} \label{GrindEQ__5_}
\left\{\begin{array}{c} {\varphi _{0} =-\tilde{u}_{0} \left(0\right),} \\ {-\varphi _{0} +\varphi _{1} =-A^{-1} \frac{d\tilde{u}_{0} \left(0\right)}{dt} ,} \\ {\varphi _{0} -2\varphi _{1} +2\varphi _{2} =-A^{-2} \frac{d^{2} \tilde{u}_{0} \left(0\right)}{dt^{2} } .} \end{array}\right.
\end{equation}
The vectors $\varphi _{0} ,\varphi _{1} ,\varphi _{2} $ are defined uniquely from the system (5) and it is clear that they belong to $D\left(A^{{\raise0.5ex\hbox{$\scriptstyle 7 $}\kern-0.1em/\kern-0.15em\lower0.25ex\hbox{$\scriptstyle 2 $}} } \right)$. As for $\left|\kappa \right|<2\lambda _{0} $  $e^{-tA} \varphi _{0} \in W_{2,\kappa }^{4} \left(R_{+} ;H\right)$, $tAe^{-tA} \varphi _{1} \in W_{2,\kappa }^{4} \left(R_{+} ;H\right)$, $t^{2} A^{2} e^{-tA} \varphi _{2} \in W_{2,\kappa }^{4} \left(R_{+} ;H\right)$, then $u(t)\in W_{2,\kappa }^{4} (R_{+} ;H)$ and satisfies the boundary-value problem (1), (2) for $A_j = 0$, $j = 1,2,3,4$.

If we take into consideration the theorem on intermediate derivatives [4, chapter 1], then the boundedness of the operator $P_{0} $ follows from the inequality

\[\left\| P_{0} u\right\| _{L_{2,\kappa } \left(R_{+} ;H\right)}^{2} =\left\| \frac{d^{4} u}{dt^{4} } +2A\frac{d^{3} u}{dt^{3} } -2A^{3} \frac{du}{dt} -A^{4} u\right\| _{L_{2,\kappa } \left(R_{+} ;H\right)}^{2} \le \]

\[4\left\| u\right\| _{W_{2,\kappa }^{4} \left(R_{+} ;H\right)}^{2} +16\left(\left\| A\frac{d^{3} u}{dt^{3} } \right\| _{L_{2,\kappa } \left(R_{+} ;H\right)}^{2} +\left\| A^{3} \frac{du}{dt} \right\| _{L_{2,\kappa } \left(R_{+} ;H\right)}^{2} \right)\le \]

\[const\left\| u\right\| _{W_{2,\kappa }^{4} \left(R_{+} ;H\right)}^{2} .\]

As a result, taking into account Banach theorem on the inverse operator, we obtain, that $P_{0} :\mathop{W_{2,\kappa }^{4} }\limits^{o} (R_{+} ;H)\to L_{2,\kappa } (R_{+} ;H)$ is the isomorphism. Theorem is proved.

\textbf{Corollary 1.} \textit{For $\left|\kappa \right|<2\lambda _{0} $ from theorem 1 it becomes clear that the norm $\left\| P_{0} u\right\| _{L_{2,\kappa } (R_{+} ;H)} $ is equivalent to the initial norm $\left\| u\right\| _{W_{2,\kappa }^{4} (R_{+} ;H)} $ in the space $\mathop{W_{2,\kappa }^{4} }\limits^{o} (R_{+} ;H)$.}

\textbf{Remark.} \textit{It is important to note that for $\kappa =\pm 2\lambda _{0} $ the operator $P_{0}$ is not invertible. Moreover, in this case $P_{0}$ is not a Fredholm operator (its image is not closed). The proof of this fact is given in [13].}

Further, we denote by $P_{1} $ the operator, acting from the space $\mathop{W_{2,\kappa }^{4} }\limits^{o} (R_{+} ;H)$ into the space $L_{2,\kappa } (R_{+} ;H)$ by the following way:

\[P_{1} u\left(t\right)\equiv \sum _{j=1}^{4}A_{j} \frac{d^{4-j} u\left(t\right)}{dt^{4-j} }  ,\, \, u\left(t\right)\in \mathop{W_{2,\kappa }^{4} }\limits^{o} \left(R_{+} ;H\right).\]

The following statement, in the proof of which we'll apply theorem on intermediate derivatives [4, chapter 1], is true.

\textbf{Lemma.} \textit{Let the operators $A_{j} A^{-j} $, $j=1,2,3,4$ be bounded in $H$. Then the operator $P_{1} $ acting from the space $\mathop{W_{2,\kappa }^{4} }\limits^{o} (R_{+} ;H)$ into the space $L_{2,\kappa } (R_{+} ;H)$ is also bounded.}

As a result we come to the main aim of this paper -- determining the sufficient conditions of regular solvability of the boundary-value problem (1), (2).

For simplifying the notes we introduce the notation $\gamma \left(\lambda \right)=1-\frac{\kappa ^{2} }{4\lambda ^{2} } $.

\textbf{Theorem 2.} \textit{Let $A=A^{*} \ge \lambda _{0} E$ $\, (\lambda _{0} >0)$, $\left|\kappa \right|<2\lambda _{0} $ and the operators $A_{j} A^{-j} $, $j=1,2,3,4$ are bounded in $H$, moreover, the inequality}

\[\sum _{j=1}^{4}c_{j} \left(\kappa \right)\left\| A_{j} A^{-j} \right\| _{H\to H}  <1\]
\textit{is true, where the numbers $c_{j} (\kappa )$, $j=1,2,3,4$ are defined by the following way:}

\[c_{1} (\kappa )=c_{3} (\kappa )=\frac{1}{2} \gamma ^{-{\raise0.5ex\hbox{$\scriptstyle 1 $}\kern-0.1em/\kern-0.15em\lower0.25ex\hbox{$\scriptstyle 2 $}} } \left(\lambda _{0} \right), c_{2} (\kappa )=\frac{1}{2\sqrt{2} } \gamma ^{-{\raise0.5ex\hbox{$\scriptstyle 1 $}\kern-0.1em/\kern-0.15em\lower0.25ex\hbox{$\scriptstyle 2 $}} } \left(\lambda _{0} \right), c_{4} \left(\kappa \right)=\gamma ^{-1} \left(\lambda _{0} \right).\]
\textit{Then the boundary-value problem (1), (2) is regularly solvable.}

\textbf{Proof.} We write the boundary-value problem (1), (2) in the form of the operator equation

\[P_{0} u\left(t\right)+P_{1} u\left(t\right)=f\left(t\right),\]
where $f(t)\in L_{2,\kappa } (R_{+} ;H)$, $u(t)\in \mathop{W_{2,\kappa }^{4} }\limits^{o} (R_{+} ;H)$. From theorem 1 it follows that the operator $P_{0} $ has the bounded inverse operator $P_{0}^{-1} $, acting from the space $L_{2,\kappa } (R_{+} ;H)$ into the space $\mathop{W_{2,\kappa }^{4} }\limits^{o} (R_{+} ;H)$. Then after the substitution $u\left(t\right)=P_{0}^{-1} z\left(t\right)$, where $z(t)\in L_{2,\kappa } (R_{+} ;H)$, we obtain the following equation in $L_{2,\kappa } (R_{+} ;H)$:

\[\left(E+P_{1} P_{0}^{-1} \right)z\left(t\right)=f\left(t\right).\]

We'll show that if the conditions of theorem are satisfied, then the norm of the operator $P_{1} P_{0}^{-1} $ is less than one. Really,

\[\left\| P_{1} P_{0}^{-1} z\right\| _{L_{2,\kappa } \left(R_{+} ;H\right)} =\left\| P_{1} u\right\| _{L_{2,\kappa } \left(R_{+} ;H\right)} \le \]

\begin{equation} \label{GrindEQ__6_}
\sum _{j=1}^{4}\left\| A_{j} \frac{d^{4-j} u}{dt^{4-j} } \right\|  _{L_{2,\kappa } \left(R_{+} ;H\right)} \le \sum _{j=1}^{4}\left\| A_{j} A^{-j} \right\| _{H\to H} \left\| A^{j} \frac{d^{4-j} u}{dt^{4-j} } \right\|  _{L_{2,\kappa } \left(R_{+} ;H\right)} .
\end{equation}
Further it is necessary to estimate the norms of intermediate derivatives operators\textit{}

\[A^{j} \frac{d^{4-j} }{dt^{4-j} } :\mathop{W_{2,\kappa }^{4} }\limits^{o} (R_{+} ;H)\to L_{2,\kappa } (R_{+} ;H), j=1,2,3,4.\]
As these operators are continuous, then their norms according to corollary 1 can be estimated with respect to $\left\| P_{0} u\right\| _{L_{2,\kappa } (R_{+} ;H)} $.

We denote by $y\left(t\right)=\left(\frac{d}{dt} +A\right)^{2} u\left(t\right)$. Then from the equation (1) for $A_j = 0$, $j = 1,2,3,4$ and boundary conditions (2) with respect to $y\left(t\right)$ we have the following boundary-value problem:\textit{}

\begin{equation} \label{GrindEQ__7_}
-\frac{d^{2} y\left(t\right)}{dt^{2} } +A^{2} y\left(t\right)=f\left(t\right),  t\in R_{+} ,
\end{equation}

\begin{equation} \label{GrindEQ__8_}
y\left(0\right)=0.
\end{equation}
After substitution $w\left(t\right)=y\left(t\right)e^{-\frac{\kappa }{2} t} $ from the problem (7), (8) we obtain

\begin{equation} \label{GrindEQ__9_}
-\left(\frac{d}{dt} +\frac{\kappa }{2} \right)^{2} w\left(t\right)+A^{2} w\left(t\right)=h\left(t\right),  t\in R_{+} ,
\end{equation}

\begin{equation} \label{GrindEQ__10_}
w\left(0\right)=0,
\end{equation}
where $w(t)\in W_{2}^{2} (R_{+} ;H)$, $h\left(t\right)=f(t)e^{-\frac{\kappa }{2} t} \in L_{2} (R_{+} ;H)$.

Multiplying both sides of the equation (9) by $A^{2} w$ as a scalar product in the space $L_{2} (R_{+} ;H)$ we have

\[\left(-\frac{d^{2} w}{dt^{2} } ,A^{2} w\right)_{L_{2} \left(R_{+} ;H\right)} +\left(-\kappa \frac{dw}{dt} ,A^{2} w\right)_{L_{2} \left(R_{+} ;H\right)} +\]

\[\left(-\frac{\kappa ^{2} }{4} w,A^{2} w\right)_{L_{2} \left(R_{+} ;H\right)} +\left(A^{2} w,A^{2} w\right)_{L_{2} \left(R_{+} ;H\right)} =\left(h,A^{2} w\right)_{L_{2} \left(R_{+} ;H\right)} .\]
Now, integrating by parts, and taking into consideration the condition (10), we obtain

\[Re\left(h,A^{2} w\right)_{L_{2} \left(R_{+} ;H\right)} =\left\| A\frac{dw}{dt} \right\| _{L_{2} \left(R_{+} ;H\right)}^{2} +\left\| A^{2} w\right\| _{L_{2} \left(R_{+} ;H\right)}^{2} -\]

\[\frac{\kappa ^{2} }{4} \left\| Aw\right\| _{L_{2} \left(R_{+} ;H\right)}^{2} \ge \left\| A\frac{dw}{dt} \right\| _{L_{2} \left(R_{+} ;H\right)}^{2} +\gamma \left(\lambda _{0} \right)\left\| A^{2} w\right\| _{L_{2} \left(R_{+} ;H\right)}^{2} \ge \]

\begin{equation} \label{GrindEQ__11_}
\gamma \left(\lambda _{0} \right)\left\| A^{2} w\right\| _{L_{2} \left(R_{+} ;H\right)}^{2} .
\end{equation}
Thus,

\[\gamma \left(\lambda _{0} \right)\left\| A^{2} w\right\| _{L_{2} \left(R_{+} ;H\right)}^{2} \le \left\| h\right\| _{L_{2} \left(R_{+} ;H\right)} \left\| A^{2} w\right\| _{L_{2} \left(R_{+} ;H\right)} ,\]
i.e.

\begin{equation} \label{GrindEQ__12_}
\left\| A^{2} w\right\| _{L_{2} \left(R_{+} ;H\right)} \le \gamma ^{-1} \left(\lambda _{0} \right)\left\| h\right\| _{L_{2} \left(R_{+} ;H\right)} .
\end{equation}
From the other side, from the inequality (11) it follows that for any $\varepsilon >0$

\[\left\| A\frac{dw}{dt} \right\| _{L_{2} \left(R_{+} ;H\right)}^{2} +\gamma \left(\lambda _{0} \right)\left\| A^{2} w\right\| _{L_{2} \left(R_{+} ;H\right)}^{2} \le \]

\[\left\| h\right\| _{L_{2} \left(R_{+} ;H\right)} \left\| A^{2} w\right\| _{L_{2} \left(R_{+} ;H\right)} \le \]

\[\frac{\varepsilon }{2} \left\| h\right\| _{L_{2} \left(R_{+} ;H\right)}^{2} +\frac{1}{2\varepsilon } \left\| A^{2} w\right\| _{L_{2} \left(R_{+} ;H\right)}^{2} .\]
We suppose in the last inequality $\varepsilon =\frac{1}{2} \gamma ^{-1} \left(\lambda _{0} \right)$ and find that

\begin{equation} \label{GrindEQ__13_}
\left\| A\frac{dw}{dt} \right\| _{L_{2} \left(R_{+} ;H\right)}^{2} \le \frac{1}{4} \gamma ^{-1} \left(\lambda _{0} \right)\left\| h\right\| _{L_{2} \left(R_{+} ;H\right)}^{2} .
\end{equation}
As $h\left(t\right)=f\left(t\right)e^{-\frac{\kappa }{2} t} $, $w\left(t\right)=y\left(t\right)e^{-\frac{\kappa }{2} t} $, then, taking into account the condition (8), from (12) and (13) we have

\[\left\| A^{2} y\right\| _{L_{2,\kappa } \left(R_{+} ;H\right)}^{2} \le \gamma ^{-2} \left(\lambda _{0} \right)\left\| f\right\| _{L_{2,\kappa } \left(R_{+} ;H\right)}^{2} ,\]

\[\left\| A\frac{dy}{dt} \right\| _{L_{2,\kappa } \left(R_{+} ;H\right)}^{2} +\frac{\kappa ^{2} }{4} \left\| Ay\right\| _{L_{2,\kappa } \left(R_{+} ;H\right)}^{2} \le \frac{1}{4} \gamma ^{-1} \left(\lambda _{0} \right)\left\| f\right\| _{L_{2} \left(R_{+} ;H\right)}^{2} .\]
From these inequalities, taking into consideration $y\left(t\right)=\left(\frac{d}{dt} +A\right)^{2} u\left(t\right)$ and the conditions (2), we obtain

\[\left\| A^{2} \frac{d^{2} u}{dt^{2} } \right\| _{L_{2,\kappa } \left(R_{+} ;H\right)}^{2} +2\left\| A^{3} \frac{du}{dt} \right\| _{L_{2,\kappa } \left(R_{+} ;H\right)}^{2} +\]

\begin{equation} \label{GrindEQ__14_}
\left\| A^{4} u\right\| _{L_{2,\kappa } \left(R_{+} ;H\right)}^{2} \le \gamma ^{-2} \left(\lambda _{0} \right)\left\| f\right\| _{L_{2,\kappa } \left(R_{+} ;H\right)}^{2} ,
\end{equation}

\[\left\| A\frac{d^{3} u}{dt^{3} } \right\| _{L_{2,\kappa } \left(R_{+} ;H\right)}^{2} +2\left\| A^{2} \frac{d^{2} u}{dt^{2} } \right\| _{L_{2,\kappa } \left(R_{+} ;H\right)}^{2} +\]

\begin{equation} \label{GrindEQ__15_}
\left\| A^{3} \frac{du}{dt} \right\| _{L_{2,\kappa } \left(R_{+} ;H\right)}^{2} \le \frac{1}{4} \gamma ^{-1} \left(\lambda _{0} \right)\left\| f\right\| _{L_{2,\kappa } \left(R_{+} ;H\right)}^{2} .
\end{equation}
As a result, from the inequalities (14) and (15) the following estimations are obtained:

\begin{equation} \label{GrindEQ__16_}
\left\| {A^j \frac{{d^{4 - j} u}}{{dt^{4 - j} }}} \right\|_{L_{2,\kappa } \left( {R_ +  ;H} \right)}  \le c_j \left( \kappa  \right)\left\| {P_0 u} \right\|_{L_{2,\kappa } \left( {R_ +  ;H} \right)} ,\,\,\,j = 1,2,3,4,
\end{equation}
where
\[c_{1} (\kappa )=c_{3} (\kappa )=\frac{1}{2} \gamma ^{-{\raise0.5ex\hbox{$\scriptstyle 1 $}\kern-0.1em/\kern-0.15em\lower0.25ex\hbox{$\scriptstyle 2 $}} } \left(\lambda _{0} \right), c_{2} (\kappa )=\frac{1}{2\sqrt{2} } \gamma ^{-{\raise0.5ex\hbox{$\scriptstyle 1 $}\kern-0.1em/\kern-0.15em\lower0.25ex\hbox{$\scriptstyle 2 $}} } \left(\lambda _{0} \right), c_{4} \left(\kappa \right)=\gamma ^{-1} \left(\lambda _{0} \right).\]

Taking into consideration the estimations (16) in the inequality (6), we have

\[\left\| P_{1} P_{0}^{-1} z\right\| _{L_{2,\kappa } \left(R_{+} ;H\right)} \le \sum _{j=1}^{4}c_{j} \left(\kappa \right)\left\| A_{j} A^{-j} \right\| _{H\to H} \left\| z\right\| _{L_{2,\kappa } \left(R_{+} ;H\right)}  .\]
Consequently,

\[\left\| P_{1} P_{0}^{-1} \right\| _{L_{2,\kappa } \left(R_{+} ;H\right)\to L_{2,\kappa } \left(R_{+} ;H\right)} \le \sum _{j=1}^{4}c_{j} \left(\kappa \right)\left\| A_{j} A^{-j} \right\| _{H\to H} <1 \]
and that's why the operator $E+P_{1} P_{0}^{-1} $ is invertible in the space $L_{2,\kappa } (R_{+} ;H)$, and it means that we can determine $u\left(t\right)$ by the formula

\[u\left(t\right)=P_{0}^{-1} \left(E+P_{1} P_{0}^{-1} \right)^{-1} f\left(t\right),\]
moreover,

\[\left\| u\right\| _{W_{2,\kappa }^{4} \left(R_{+} ;H\right)} \le \left\| P_{0}^{-1} \right\| _{L_{2,\kappa } \left(R_{+} ;H\right)\to W_{2,\kappa }^{4} \left(R_{+} ;H\right)} \times \]

\[\left\| \left(E+P_{1} P_{0}^{-1} \right)^{-1} \right\| _{L_{2,\kappa } \left(R_{+} ;H\right)\to L_{2,\kappa } \left(R_{+} ;H\right)} \left\| f\right\| _{L_{2,\kappa } \left(R_{+} ;H\right)} \le \]

\[const\left\| f\right\| _{L_{2,\kappa } \left(R_{+} ;H\right)} .\]
Theorem is proved.

\textbf{Corollary 2. }\textit{Let $\kappa =0$ and the inequality}

\[\frac{1}{2} \left\| A_{1} A^{-1} \right\| _{H\to H} +\frac{1}{2\sqrt{2} } \left\| A_{2} A^{-2} \right\| _{H\to H} +\frac{1}{2} \left\| A_{3} A^{-3} \right\| _{H\to H} +\]

\[\left\| A_{4} A^{-4} \right\| _{H\to H} <1\]
\textit{is true. Then the operator $P_{0} +P_{1} $ is the isomorphism from the space $\mathop{W_{2}^{4} }\limits^{o} (R_{+} ;H)$ onto the space $L_{2} (R_{+} ;H)$.}

\bigskip
\newpage

\begin{center}{\large \bf  REFERENCES}
\end{center}

[1] Gorbachuk, M.L. and Gorbachuk, V.I. On Well-Posed Solvability in Some Classes of Entire Functions of the Cauchy Problem for Differential Equations in a Banach Space // \textit{Methods of Functional Analysis and Topology, 2005, Vol. 11, No. 2, pp. 113-125.}

[2] Favini, A. and Yakubov, Ya. Regular Boundary Value Problems for Elliptic Differential-Operator Equations of the Fourth Order in UMD Banach Spaces // \textit{Scientiae Math. Japonicae, 2009, Vol. 70, No. 2, pp. 183-204.}

[3] Aliev, A R. and Mirzoev, S.S. On Boundary Value Problem Solvability Theory for a Class of High-Order Operator-Differential Equations // \textit{Functional Analysis and Its Applications}, \textit{2010}, \textit{Vol}. \textit{44}, \textit{No}. \textit{3}, \textit{pp}. \textit{209-211}. (\textit{published in Funktsional'nyi Analiz i Ego Prilozheniya}, \textit{2010, Vol. 44, No. 3, pp. 63-65.})

[4] Lions, J.-L. and Magenes, E. \textit{Non-Homogeneous Boundary Value Problems and Applications}, Paris: Dunod, 1968; Moscow: Mir, 1971; Berlin: Springer, 1972.

[5] Mirzoev, S.S. On Solvability of Boundary-Value Problems for Operator-Differential Equations of the Second Order in Spaces with Weight / in \textit{Lineinye operatory i ikh prilozheniya }(\textit{Linear Operators and Their Applications}), \textit{Baku, 1989, pp. 46-49.} (Russian)

[6] Aliev, A.R. Solubility of Boundary-Value Problems for a Class of Third-Order Operator-Differential Equations in a Weighted Space // \textit{Russian Math. Surveys, 2005, Vol. 60, No. 4, pp. 791-793. } (\textit{published in Uspekhi Mat. Nauk, 2005, Vol. 60, No. 4 (364), pp. 215-216.})

[7] Aliyev, A.R. On   Solvability of a Boundary-Value Problem for Elliptic Type Operator-Differential Equations with Discontinuous Coefficients in the Weight Space // \textit{Proceedings of Institute of Mathematics and Mechanics of NAS of Azerbaijan, 2006, Vol. 24 (32), pp. 17-28.}

[8] Aliev, A.R. Solvability of a Class of Boundary Value Problems\textit{ }for Second-Order Operator-Differential Equations\textit{ }with a Discontinuous Coefficient in a Weighted Space // \textit{Differential Equations, 2007, Vol. 43, No. 10, pp. 1459-1463. }(\textit{published in Differentsial'nye Uravneniya, 2007, Vol. 43, No. 10, pp. 1423-1426.})

[9] Aliev, A.R. and Gasymov, A.A.  On the Correct Solvability of the Boundary-Value Problem for One Class Operator-Differential Equations of the Fourth Order with Complex Characteristics // \textit{Boundary Value Problems, 2009, Vol. 2009, Article ID 710386, 20 pages.}

[10] Gumbataliev, R.Z. Normal Solvability of Boundary Value Problems for a Class of Fourth-Order Operator-Differential Equations in a Weighted Space // \textit{Differential Equations, 2010, Vol. 46, No. 5, pp. 681-689.} (\textit{published in Differentsial'nye Uravneniya, 2010, Vol. 46, No. 5, pp. 678-686.})

[11] Dubinskii, Yu.A. On Some Differential-Operator Equations of Arbitrary Order // \textit{Mathematics of the USSR-Sbornik, 1973, Vol. 19, No. 1, pp. 1-21.} (\textit{published in Mat. Sbornik,  1973, Vol. 90 (132), No. 1, pp. 3-22.})

[12] Shkalikov, A.A. Some Problems in the Theory of Polynomial Operator Pencils // \textit{Russian Mathematical Surveys, 1983, Vol. 38, No. 3, pp. 151-152.} (\textit{published in Uspekhi Mat. Nauk, 1983, Vol. 38, No. 3 (231), pp. 189–190.})

[13] Shkalikov, A.A. Elliptic Equations in Hilbert Space and Associated Spectral Problems // \textit{J. Soviet Math., 1990, Vol. 4, pp. 2399-2467.} (\textit{published in Trudy Sem. Petrovsk., 1989, No. 14, pp. 140-224.})

[14] Teters, G.A. \textit{Complex Loading and Stability of the Covers from Polymeric Materials}, Latvia, Riga: Zinatne Press, 1969. (Russian)

\bigskip

\newpage

\begin{center}{\large \bf On the Boundary-Value Problem for One Class of Differential Equations of the Fourth Order with the Operator Coefficients}
\end{center}

\begin{center}{\bf A.R.~Aliev}
\\ {\it Institute~of~Mathematics~and~Mechanics~of~NAS~of~Azerbaijan}
\\ {\it 9, F. Agayev Str., AZ1141, Baku, Azerbaijan,}
\\ {\it Baku~State~University}
\\ {\it 23, Z. Khalilov Str., AZ1148, Baku, Azerbaijan}
\\ {\it E-mail: alievaraz@yahoo.com}
\end{center}

\begin{quote}
\small {\bf Abstract.} The boundary-value problem on semi-axis for one class operator-differenti-al equations of the fourth order, the main part of which has the multiple characteristic is investigated in this paper in Sobolev type weighted space. Correctness and unique solvability of the boundary-value problem is proved, and the solvability conditions are expressed in terms of the operator coefficients of the equation. Estimations of the norms of the operators of intermediate derivatives, closely connected with the solvability conditions, have been carried out. The connection between the exponent of the weight and the lower border of the spectrum of the main operator, participating in the equation, is determined in the results of the paper.
\end{quote}
\medskip

\begin{quote}
{\small {\bf 2000 Mathematics Subject Classification:} 34G10, 34K10, 35J40, 47D03.}

{\small {\bf Key words and phrases:} boundary-value problem, selfadjoint operator, Hilbert space, weighted space, operator-differential equation, complex characteristic, regular solvability, isomorphism, intermediate derivatives.}
\end{quote}
\medskip

\bigskip

\end{document}